\documentclass[leqno]{article}
\usepackage{amssymb}
\usepackage{amsmath, enumerate, vmargin}

\makeatletter
\renewcommand\section{\@startsection {section}{1}{\z@}%
 {-3.5ex \@plus -1ex \@minus -.2ex}%
 {2.3ex \@plus.2ex}%
 {\center \normalfont\large\bfseries}}
\makeatother

\setmarginsrb{3cm}{3cm}{3cm}{3cm}{1mm}{6mm}{0mm}{10mm}

\newtheorem{thm}{Theorem}[section]

\newtheorem{cor}[thm]{Corollary}
\newtheorem{lem}[thm]{Lemma}
\newtheorem{defi}[thm]{Definition}
\newtheorem{remark}[thm]{Remark}
\newtheorem{example}[thm]{Example}
\newtheorem{pb}[thm]{Problem}

\newenvironment{rk}{\begin{remark}\rm}{\end{remark}}

\newcommand{\real}{{\mathbb R}}
\newcommand{\nat}{{\mathbb N}}

\newcommand{\com}{{\mathbb C}}
\newcommand{\un}{1\mkern -4mu{\textrm l}}

\newcommand{\E}{{\mathbb E}}
\renewcommand{\H}{{\mathcal H}}
\newcommand{\M}{{\mathsf M}}
\newcommand{\Tr}{{\rm Tr}}

\renewcommand{\a}{\alpha}
\renewcommand{\b}{\beta}

\renewcommand{\d}{\delta}
\newcommand{\g}{\gamma}
\newcommand{\e}{\varepsilon}
\renewcommand{\O}{\Omega}
\renewcommand{\o}{\omega}
\renewcommand{\l}{\lambda}
\newcommand{\La}{\Lambda}
\renewcommand{\e}{\varepsilon}
\newcommand{\f}{\varphi}

\newcommand{\ot}{\otimes}

\renewcommand{\t}{\tau}
\newcommand{\8}{\infty}
\newcommand{\el}{\ell}
\newcommand{\la}{\langle}
\newcommand{\ra}{\rangle}
\newcommand{\wt}{\widetilde}

\newcommand{\n}{\noindent}
\newcommand{\pf}{\noindent{\it Proof.~~}}
\newcommand{\cqd}{\hfill$\Box$}
\newcommand{\be}{\begin{eqnarray*}}
\newcommand{\ee}{\end{eqnarray*}}
\newcommand{\beq}{\begin{equation}}
\newcommand{\eeq}{\end{equation}}

\numberwithin{equation}{section}

\begin{document}


\title{The little Grothendieck theorem and Khintchine inequalities \\
for symmetric spaces of measurable operators}
\author{Fran\c{c}oise Lust-Piquard and Quanhua Xu}

\date{}

\maketitle


\begin{abstract}
 We prove the little Grothendieck theorem for any 2-convex
noncommutative symmetric space. Let $\M$ be a von Neumann algebra
equipped with a normal faithful semifinite trace $\t$, and let $E$
be an r.i. space on $(0,\;\8)$. Let $E(\M)$ be the associated
symmetric space of measurable operators. Then to any bounded
linear map $T$ from $E(\M)$ into a Hilbert space $\mathcal H$
corresponds a positive norm one functional $f\in E_{(2)}(\M)^*$
such that
 $$\forall\; x\in E(\M)\quad \|T(x)\|^2\le K^2\,\|T\|^2 f(x^*x+xx^*),$$
where $E_{(2)}$ denotes the 2-concavification of $E$ and  $K$ is a
universal constant. As a consequence we obtain the noncommutative
Khintchine inequalities for $E(\M)$ when  $E$ is either 2-concave
or 2-convex and $q$-concave for some $q<\8$. We apply these
results to the study of Schur multipliers from a 2-convex unitary
ideal into a 2-concave one.

\end{abstract}



\makeatletter
\renewcommand{\@makefntext}[1]{#1}
\makeatother \footnotetext{\noindent
 F.P.: D{\'e}partement de Math\'{e}matiques, Universit\'{e} de
 Cergy-Pontoise,
  95302, Cergy-Pontoise, cedex - France\\
  francoise.piquard@math.u-cergy.fr\\
Q.X.:  Laboratoire de Math\'{e}matiques, Universit\'{e} de
Franche-Comt\'{e},
25030 Besan\c con, cedex - France\\
 qx@math.univ-fcomte.fr
 }


 \section{Introduction}


Let $C(\O)$ denote the space of continuous functions on a compact
topological space $\O$, equipped with the uniform norm. The
classical little Grothendieck theorem asserts that for any bounded
linear map $T$ from $C(\O)$ into a Hilbert space $\H$ there exists
a probability measure $\mu$ on $\O$ such that
 $$\forall\; x\in C(\O)\quad \|T(x)\|^2\le
 K^2\,\|T\|^2\int_\O|x|^2\,d\mu\,,$$
where $K$ is an absolute positive constant. This result was
extended by Maurey \cite{mau-ast} to maps defined on any 2-convex
Banach lattice $\La$. Namely, if $T:\La\to \H$ is bounded, then
there exists a positive norm one functional $f\in (\La_{(2)})^*$
such that
 $$\forall\; x\in\La\quad \|T(x)\|^2\le
 K^2\,\|T\|^2f(|x|^2)\,.$$
Here $\La_{(2)}$ denotes the 2-concavification of $\La$. The
reader is referred to \cite{LT-II} for all notions on Banach
lattices used in this paper.

On the other hand, the noncommutative analogue of the little
Grothendieck theorem was obtained by Pisier \cite{pis-groC*} (see
also \cite{pis-fact}). More precisely, let $A$ be a C*-algebra,
and let $T:A\to \H$ be a bounded linear map. Then there exists a
state $f$ on $A$ such that
 $$\forall\; x\in A\quad \|T(x)\|^2\le
 K^2\,\|T\|^2f(x^*x+xx^*)\,.$$
In the spirit of Pisier's theorem, the first named author of the
present paper extended in \cite{lust-gro} Maurey's inequality to
unitary ideals of operators on a Hilbert space, and more
generally, to symmetric spaces of measurable operators, provided
that the underlying r.i. spaces are 2-convex and satisfy an
additional condition (see the discussion following Theorem
\ref{gro thm} below for more details). It was conjectured in
\cite{lust-gro} that this additional condition should be
irrelevant.

\medskip

The main objective of this paper is to remove the additional
condition mentioned above from the main result of \cite{lust-gro},
so we obtain the full noncommutative analogue of Maurey's
inequality. On the other hand, the arguments of \cite{lust-gro}
are rather lengthy, and unfortunately, contain some obscure points
about polar decomposition (see \cite[Lemma IV.5]{lust-gro}). Our
proof of Theorem \ref{gro thm} below is simpler and more readable.
To state our main result we need to introduce symmetric spaces of
measurable operators.

\medskip

Let $\mathsf{M}$  be a von Neumann algebra, equipped with a
semifinite normal faithful trace $\t$, and  let $L_0(\M, \t)$, or,
simply $L_0(\M)$ denote the topological $*$-algebra of all
operators which are measurable with respect to $(\M,\t)$. The
topology of $L_0(\M)$ is determined by convergence in measure. For
$x\in L_0(\M)$ and $t>0$, $\mu_t(x)$ denotes the t-th generalized
singular number of $x$. The function $t\mapsto \mu_t(x)$ is called
the generalized singular number function and is denoted by
$\mu(x)$. Recall that $\mu(x)$ is nonincreasing and
$\mu(x)=\mu(x^*)=\mu(|x|)$, where $|x|=(x^*x)^{1/2}$ is the
absolute value of $x$. The reader is referred to \cite{fk} for
more details on generalized singular numbers.

Let $E$ be an r.i. space on $(0,\8)$ in the sense of \cite{LT-II}.
The symmetric space $E(\M,\t)$ of measurable operators associated
with $\M$ and $E$ is defined as the space of all measurable
operators $x\in L_0(\M)$ such that $\mu(x)\in E$. $E(\M,\t)$ is a
Banach space equipped with the norm
$\|x\|_{E(\M,\t)}=\|\mu(x)\|_E$. $E(\M,\t)$ is often denoted
simply by $E(\M)$. The spaces $E(\M)$ are the so-called
noncommutative symmetric spaces, studied in detail for the first
time by Ovchinnikov \cite{ovchi}.
 Note that if $\M=B(\el_2)$
and $\t$ is the usual trace on $B(\el_2)$,  $E(\M)$ is a unitary
ideal of operators on $\el_2$. On the other hand, if $\t$ is
finite, $E$ can be taken to be an r.i. space on $[0,\:\t(1)]$.
Recall that if $E=L_p(0,\;\8)$,  $E(\M)=L_p(\M)$, the
noncommutative $L_p$-space associated with $(\M,\t)$.

For $r>1$, $E^{(r)}$  and $E_{(r)}$ denote the $r$-convexification
and $r$-concavification of $E$, respectively. Recall that if $E$
is a $p$-convex and $q$-concave r.i. space,  $E^{(r)}$ is a
$pr$-convex and $qr$-concave r.i. space. If in addition $p\ge r$
and the $p$-convexity constant of $E$ is equal to 1, then
$E_{(r)}$ is a $p/r$-convex and $q/r$-concave r.i. space. In
particular, if $E$ is 2-convex with constant 1, $E_{(2)}$ is an
r.i. space. Recall that if $E$ is $p$-convex and $q$-concave, $E$
can be renormed into an r.i. space which is  $p$-convex and
$q$-concave with constant 1.

\medskip

The following is our main result. In the remainder of the paper,
unless explicitly stated otherwise, $\M$ will denote a von Neumann
algebra equipped with a normal faithful semifinite trace $\t$, and
$E$ will be an r.i. space on $(0, \;\8)$. $K$ will denote a
universal positive constant, which may change from line to line.

\begin{thm}\label{gro thm}
 Assume that $E$ is  $2$-convex with constant $1$.
Let $\H$ be a Hilbert space. Then for any bounded linear map $T:
E(\M)\to \H$ there exists a positive norm one functional $f\in
E_{(2)}(\M)^*$ such that
 $$\forall\; x\in E(\M)\quad \|T(x)\|^2\le
 K^2\,\|T\|^2f(x^*x+xx^*)\,.$$
 \end{thm}

This theorem is stated in \cite{lust-gro} with the stronger
assumption that $E$ is $p$-convex with $p>2$. For unitary ideals
(i.e. when $\M=B(\el_2)$ equipped with the usual trace), the
$p$-convexity assumption is weakened to 2-convexity plus an
additional mild condition.

We should also emphasize the universality of the constant $K$ in
Theorem \ref{gro thm}, which is of independent interest.  In some
special cases, it is much easier to prove the little Grothendieck
inequality with a constant depending on the space $E$ in
consideration. This is, for instance, the case for $E=L_p(0,\;\8)$
with $2\le p<\8$ (see \cite[Theorem 6.6]{px-survey}).

\medskip

The proof of Theorem \ref{gro thm}  will be given in the next
section. It depends on two other equivalent statements. One of
them is the (difficult) lower estimate in the noncommutative
Khintchine inequalities for the dual space $E(\M)^*$ of $E(\M)$,
which is important for its own right. To state this equivalence it
is more convenient to work with the noncommutative symmetric space
$E'(\M)$ instead of $E(\M)^*$, where $E'$ denotes the K\"{o}the
dual of $E$, which is the subspace of $E^*$ consisting of all
integrals. Let us recall the well-known relations between $E^*$
and $E'$. $E'$ is a norming subspace of $E^*$. If $E$ is order
continuous, $E^*=E'$. On the other hand, if $E$ is maximal (i.e.
$E=E''$) and is $p$-convex with $p>1$, then $E={E'}^*$. Indeed,
the $p$-convexity of $E$ implies that the restriction of $E$ to
$[0,\;1]$ is not order isomorphic to $L_1[0,\;1]$. Thus by
Proposition 2.a.3 and the remark following it in \cite{LT-II},
$E=F^*$, where $F$ is the closure of simple functions in $E'$.
However, since $E'$ is $p'$-concave ($p'$ denoting the conjugate
index of $p$),  $E'$ is order continuous. It follows that $E'=F$.
Also observe that if $E$ is an r.i. space on $[0, 1]$ and is not
order continuous, then $E$ is maximal. Indeed, if $E$ is minimal
and non separable, $E$ is order isomorphic to $L_\8(0,1)$, so
maximal. Consequently, every r.i. space on $[0, 1]$ either is
order continuous or has the Fatou property.

Now if $E$ is order continuous, $E(\M)^*=E^*(\M)$. This is
\cite[Lemma 1]{xu-RN} if $\t(1)<\8$ and is stated in \cite[p.
745]{ddp} for the general case. Let us also note that the latter
case easily follows from the former by a standard approximation
argument using the semifiniteness of $\t$.

\medskip

We are ready to state  the equivalence theorem of \cite{LPP}.
$(\e_k)$ denotes a Rademacher sequence on a probability space, and
$\E$ is the corresponding expectation.

\begin{thm}\label{equivalence}
 Assume that $E$ is  $2$-convex with constant 1.
Then the following statements are equivalent
 \begin{enumerate}[\rm i)]
 \item
There exists  a positive constant $C_1$ such that to any bounded
map $T$ from $E(\M)$ into a Hilbert space $\H$ corresponds a
positive norm one functional $f\in E_{(2)}(\M)^*$ satisfying
 $$
 \forall\; x\in E(\M)\quad
 \|T(x)\|\le C_1\|T\|\big(f(x^*x+xx^*)\big)^{1/2}\,.\leqno (\rm G)
 $$
 \item There exists  a positive constant $C_1$ such that
for any bounded map $T: E(\M)\to \H$ and any finite sequence
$(x_k)\subset E(\M)$
 $$\big(\sum_k\|T(x_k)\|^2\big)^{1/2}\le C_1\|T\| \big\|\big(
 \sum_k x_k^*x_k + x_kx_k^*\big)^{1/2}\big\|\,.\leqno(\rm C)$$
 \item There exists a positive constant $C_2$ such
that for any finite sequence $(x_k)\subset E'(\M)$
 $$\inf\big\{\big\|\big(\sum_ka_k^*a_k\big)^{1/2}\big\|+
 \big\|\big(\sum_kb_kb_k^*\big)^{1/2}\big\|\big\}\le C_2
 \big(\E \big\|\sum_k\e_k\,x_k\big\|^2\big)^{1/2},\leqno(\rm K)$$
where the infimum runs over all decompositions $x_k=a_k+b_k$ in
$E'(\M)$.
 \end{enumerate}
Moreover, the constants $C_1$ and $C_2$ above satisfy the
relations: $C_1\le C_2\le KC_1$.
 \end{thm}

By standard arguments we obtain the following general
noncommutative Khintchine inequalities for symmetric spaces of
measurable operators. They generalize the Khintchine inequalities
for noncommutative $L_p$-spaces in \cite{lust-khin} and
\cite{LPP}.

 \begin{thm}\label{khintchine}
  \begin{enumerate}[\rm i)]
  \item If $E$ is $2$-concave with constant 1,
 then for every finite sequence $(x_k)\subset E(\M)$
  \be
 \big(\E \big\|\sum_k\e_k\,x_k\big\|^2\big)^{1/2}
 \le\inf\big\{\big\|\big(\sum_ka_k^*a_k\big)^{1/2}\big\|+
 \big\|\big(\sum_kb_kb_k^*\big)^{1/2}\big\|\big\}
 \le K\big(\E \big\|\sum_k\e_k\,x_k\big\|^2\big)^{1/2},
 \ee
where the infimum runs over all decompositions $x_k=a_k+b_k$ in
$E(\M)$.
 \item If $E$ is $2$-convex and $q$-concave with constant 1
for some $q<\8$, then for every finite sequence $(x_k)\subset
E(\M)$
  \be
 K_q^{-1}\big(\E \big\|\sum_k\e_k\,x_k\big\|^2\big)^{1/2}
 \le \max\big\{\big\|\big(\sum_kx_k^*x_k\big)^{1/2}\big\|,\
 \big\|\big(\sum_kx_kx_k^*\big)^{1/2}\big\|\big\}
 \le \big(\E \big\|\sum_k\e_k\,x_k\big\|^2\big)^{1/2}\,,
 \ee
where $K_q$ depends only on $q$. Moreover, $K_q\le K q$.
 \end{enumerate}
 \end{thm}

\pf The second inequality of i) follows from Theorem \ref{gro thm}
and Theorem \ref{equivalence}. The first one is obtained by using
the 2-concavity of $E$ as in \cite{LPP}. In the same way, the
second inequality of ii) is a consequence of the 2-convexity of
$E$. Thus it remains to prove the first one of ii). This is done
via duality by using the second inequality of i). To this end, we
need the K-convexity of $E(\M)$ (cf. e.g. \cite{pis-convbody} for
the definition of K-convexity) . Under the assumption of ii), by
\cite{xu-sym}, $E(\M)$ is of type 2 with a constant $T_q$
depending only on $q$, so $E(\M)$ is K-convex. Alternately, we can
also use \cite[Theorem 7.11]{pis-convbody}. Indeed, by
\cite{pis-comint}, there exists an r.i. space $E_0$ such that
$E=\big(E_0,\;L_2(0,\;\8)\big)_\theta$, where $\theta=2/q$. Then
$E(\M)=\big(E_0(\M),\;L_2(\M)\big)_\theta$. Thus it follows that
$E(\M)$ is K-convex with constant majorized by $K'q$ for some
universal constant $K'$. Therefore, using the second inequality of
ii) and duality, we deduce the first inequality of ii) with
$K_q\le KK'q$. \cqd

\medskip

Note that the $q$-concavity condition in Theorem \ref{khintchine},
ii) is necessary. Indeed, under the 2-convexity assumption of $E$,
the first inequality of ii) implies that $E$ is of type 2, and so
is of finite concavity. On the other hand, if $E=L_q(0, \8)$ with
$2\le q<\8$, the optimal order of the constant $K_q$ above is
${\rm O}(\sqrt q)$. We do not know whether this is true in the
general case.

\medskip

We end this section with some open problems. The first one
concerns the noncommutative Khintchine inequalities. Theorem
\ref{khintchine} gives a deterministic characterization of the
expression $\E \big\|\sum_k\e_k\,x_k\big\|$ only when $E$
satisfies one of the two conditions there. Recall that if $E$ is a
$q$-concave Banach lattice for some $q<\8$, then for any finite
sequence $(x_k)\subset E$
 $$\E \big\|\sum_k\e_k\,x_k\big\|\approx
 \big\|\big(\sum_k|x_k|^2\big)^{1/2}\big\|$$
with relevant constants depending only on $q$ and the
$q$-concavity constant of $E$. At the time of this writing, we do
not know how to characterize deterministically $\E
\big\|\sum_k\e_k\,x_k\big\|_{E(\M)}$ for a general $E$.

\begin{pb}
 Let $E$ be a $q$-concave r.i. space with $q<\8$. Find a
deterministic characterization of $\E \big\|\sum_k\e_k\,x_k\big\|$
for any finite sequence $(x_k)\subset E(\M)$.
 \end{pb}

The second  problem is on the big Grothendieck theorem in the
setting of this paper.

\begin{pb}
 Let $E$ and $F$ be two  $2$-convex r.i. spaces with constant $1$.
Let $u: E(\M)\times F(\M)\to\com$ be a bounded bilinear form. Do
there exist two positive norm one functionals $f\in E_{(2)}(\M)^*$
and $g\in F_{(2)}(\M)^*$ such that
 $$\forall\; x\in E(\M),\; \forall\;y\in F(\M)\quad
 |u(x, y)|\le K\|u\| \big[f(x^*x+xx^*)\big]^{1/2}\,
 \big[g(y^*y+yy^*)\big]^{1/2}\,?$$
 \end{pb}

We can state the following more general problem.

\begin{pb}
 Let $E$  be a  $2$-convex r.i. space with constant $1$ and $Y$ a
Banach space of cotype $2$. Let $T: E(\M)\to Y$ be a bounded
linear map. Do there exist a positive norm one functional $f\in
E_{(2)}(\M)^*$  and a positive constant $C$ $($depending only on
the cotype $2$ constant of $Y)$ such that
 $$\forall\; x\in E(\M)\quad
 \|T(x)\|^2\le C^2\|T\|^2 f(x^*x+xx^*)\,?$$
This can be reformulated as follows. Does there exist  a positive
constant $C$ $($depending only on the cotype $2$ constant of $Y)$
such that
 $$\forall\; x_k\in E(\M)\quad
 \big(\sum_k \|T(x_k)\|^2\big)^{1/2}\le C\|T\|
 \big\|\big(\sum_kx_k^*x_k+x_kx_k^*\big)^{1/2}\big\|\,?$$
 \end{pb}

In the case of $E=L_\8(0,\8)$ (i.e. $E(\M)=\M$) the previous
problem has a positive solution. In this case $\M$ can be replaced
by any C*-algebra (see \cite{pis-mathann}).


 \section{Proof of Theorem \ref{gro thm}}

This section  is devoted to the proof of Theorem \ref{gro thm}. We
require two auxiliary lemmas.

 \begin{lem}\label{approximation}
 Assume that $\t(1)=1$ and $E=F^*$ for an order continuous r.i.
space $F$ on $[0,\;1]$.
 \begin{enumerate}[\rm i)]
 \item $E(\M)=F(\M)^*$.
 \item For any $T:  E(\M)\to \H$  with $\|T\|\le1$
there exists a net $(T_i)$ of w*-continuous finite rank maps from
$E(\M)$ into $\H$ such that $\|T_i\|\le 1$ and $T_i\to T$ strongly
$($i.e. in the point-norm topology$)$.
 \end{enumerate}
 \end{lem}

\pf i)  is \cite[Lemma 1]{xu-RN}. To show ii) we use standard
duality for Banach space tensor product. We have
 $$B(E(\M),\; \H)=\big(E(\M){\mathop\ot^\wedge} \H^*\big)^*\,,$$
where $\displaystyle\mathop\ot^\wedge$ denotes the projective
tensor product for Banach spaces, and where the duality is
determined as follows. For $T\in B(E(\M),\; \H)$ and $x\ot \xi\in
E(\M)\ot \H^*$
 $$\la T,\;x\ot \xi\ra=\la\xi,\; T(x)\ra.$$
On the other hand,
 $$E(\M){\mathop\ot^\wedge} \H^*
 =\big(F(\M){\mathop\ot^\vee} \H\big)^*\,,$$
where $\displaystyle\mathop\ot^\vee$ denotes the injective tensor
product for Banach spaces. It follows that $B(E(\M),\; \H)$ is the
bidual of $F(\M){\displaystyle\mathop\ot^\vee} \H$. Therefore, the
unit ball of $B(E(\M),\; \H)$ is the w*-closure of that of
$F(\M){\displaystyle\mathop\ot^\vee} \H$. Recall that
$F(\M){\displaystyle\mathop\ot^\vee} \H$ can be identified as the
norm closure in $B(E(\M),\; \H)$ of w*-continuous finite rank maps
(i.e. those associated with vectors in the algebraic tensor
product $F(\M)\ot \H$). Now let $T: E(\M)\to \H$ with $\|T\|\le
1$. Then we deduce a net $(T_i)$ such that each $T_i$ is a
w*-continuous finite rank map from $E(\M)$ to $\H$, $\|T_i\|\le 1$
and $T_i\to T$ in the w*-topology. Thus $T_i(x)\to T(x)$ weakly in
$\H$ for any $x\in E(\M)$. Therefore, an appropriate net of convex
combinations of the $T_i$'s converges to $T$ strongly. \cqd

\begin{lem}\label{density}
 Assume $\t(1)=1$ and $E$ is an order continuous
r.i. space on $[0,1]$. Let $a\in E(\M)$ be a positive operator and
$e=s(a)$ the support projection of $a$. Then $\{ha+ah +e^\perp
he^\perp:\; h\in \M_h\}$ is dense in $E(\M)_h$, where
$e^\perp=1-e$ and $E(\M)_h$ denotes the selfadjoint part of
$E(\M)$, i.e. $E(\M)_h=\{x\in E(\M)\;:\; x^*=x\}$.
 \end{lem}

\pf The order continuity of $E$ implies that $E^*=E'$ is again an
r.i. space on $[0,1]$. Thus by Lemma \ref{approximation} i),
$E(\M)^*=E^*(\M)=E'(\M)$. Let $y\in E^*(\M)$ be such that
 $$\forall\; h\in \M_h\quad
 \t\big(y(ha+ah +e^\perp he^\perp)\big)=0.$$
Then
 \beq\label{density1}
 \forall\; x\in \M\quad
 \t\big((ya+ay)x +e^\perp ye^\perp x)\big)=0.
 \eeq
This implies in particular
 $$\forall\; x\in \M\quad\t\big((e^\perp y e^\perp)x\big)=0;$$
whence
 \beq\label{density2}
 e^\perp y e^\perp=0.
 \eeq
Thus by (\ref{density1}) we deduce that $ya=-ay$, so $ya^2=a^2y$.
Therefore, $y$ commutes with all polynomials in $a^2$,  thus by
functional calculus, with $(a^2)^{1/2}=a$ too. It follows that
$ay=ya=0$; whence $ey=ye=0$. Combining this with (\ref{density2}),
we get $y=0$, which implies the desired density.\cqd

\medskip

\n{\it Proof of Theorem} \ref{gro thm}. We will prove one of the
three equivalent statements in Theorem \ref{equivalence} according
to different cases. We start the proof by reducing $\t$ to a
finite trace. To this end we consider the noncommutative
Khintchine inequality in Theorem \ref{equivalence}. Note that $E'$
is 2-concave for $E$ is 2-convex. Thus $E'$ is order continuous.
On the other hand, by the semifiniteness of $\t$, we have an
increasing net $(e_i)$ of projections in $\M$ such that
$\t(e_i)<\8$ for each $i$ and $e_i\to1$ strongly. Then the order
continuity of $E'$ implies that $e_ixe_i\to x$ in $E'(\M)$ for
every $x\in E'(\M)$ (see \cite[Lemma 4.5]{xu-camb}). Therefore, we
need only to prove inequality (K) for all $x_k$ such that
$x_k=e_ix_ke_i$ and for every fixed $i$. Namely, it suffices to
show (K) for $E'(e_i\M e_i)$. Now the restriction of $\t$ to
$e_i\M e_i$ is finite. Thus we are reduced to the finite trace
case.

In the sequel, $\t$ is  a normal faithful finite trace on $\M$, so
by normalization, we can further assume $\t(1)=1$. Accordingly,
$E$ can be taken to be a 2-convex r.i. space on $[0,1]$. The
remainder of the proof is divided into several cases.

\medskip

\n {\bf Case 1:} {\it $E$ is  $p$-convex and $q$-concave with
constant $1$ for some $p>2$ and $q<\8$}.
 This is the main part of the whole proof.
We will prove the little Grothendieck theorem for $E(\M)$. The
pattern of the following argument is modelled on Haagerup's proof
of the little Grothendieck theorem for C*-algebras (see
\cite{haag-gro}). It is clear that it suffices to do this for
every finite dimensional Hilbert space $\H$. So in this part $\H$
is assumed finite dimensional. Fix a map $T:\; E(\M)\to \H$ such
that $\|T\|=1$. The $p$-convexity and $q$-concavity of $E$ implies
that $E(\M)$ is uniformly convex and uniformly smooth by virtue of
\cite{xu-sym}. In particular, $E(\M)$ is reflexive. Then $T$ is
weakly continuous, so the weak compactness of the unit ball of
$E(\M)$ implies that $T$ attains its norm (recalling that $\dim
\H<\8$). Thus there exists $a\in E(\M)$ such that $\|a\|=1$ and
$\|T(a)\|=1$. We consider two subcases according to $a\ge0$ or
not.

\medskip\n{\bf Subcase 1: $a\ge0$.}
  Let $h\in \M$ be a
selfadjoint operator. Then $e^{ith}$ is  unitary for any
$t\in\real$. Consequently,
 $$e^{ith}a\,e^{ith}\in E(\M)\quad \mbox{and}\quad
 \|e^{ith}a\,e^{ith}\|=1.$$
Writing
 $$e^{ith}a\,e^{ith}=a-t^2b +it(ha+ah)+ {\rm o}(t^2),$$
where $b=(h^2a + ah^2)/2 + hah$, we have
 \be
 \|a-t^2b +it(ha+ah)\|\le 1 +{\rm o}(t^2).
 \ee
By the selfadjointness of $h$,
 \be
 \|a-t^2b +it(ha+ah)\|=\|a-t^2b -it(ha+ah)\|.
 \ee
Thus
 $$\E\|a-t^2b +it\e (ha+ah)\|^2
 = \|a-t^2b +it(ha+ah)\|^2,$$
where $\e$ is a Rademacher function and $\E$ the corresponding
expectation. Then we deduce that
 \be
 \|T(a-t^2b)\|^2 +t^2\|T(ha+ah)\|^2
 &=& \E\|T(a-t^2b) +it\e T(ha+ah)\|^2\\
 &\le&\E\|a-t^2b +it\e (ha+ah)\|^2
 \le 1 +{\rm o}(t^2).
 \ee
Therefore,
 $$t^2\|T(ha+ah)\|^2\le 2t^2\,{\rm Re}\la T(a), \;T(b)\ra
 +{\rm o}(t^2);$$
whence
 \beq\label{gro thm1}
 \|T(ha+ah)\|^2\le2\, {\rm Re}\la T(a), \;T(b)\ra.
 \eeq

Let $\f=T^*T(a)$. (More rigorously, $\f=\overline{T^*}\,T(a)$ with
$\overline{T^*}: H=\overline{H^*}\to \overline{E(\M)^*}$, where
$\overline X$ denotes the complex conjugate of a Banach space
$X$.)  Then $\f\in E(\M)^*=E^*(\M)$ and $\|\f\|\le 1$. On the
other hand, $\f(a)=1$. Consequently, $\|\f\|=1$ and $\f$ is a
supporting functional of $a$, which is unique by virtue of the
smoothness of $E(\M)$. $\f$ must be positive and $s(\f)\le e$,
where  $e=s(a)$ is the support projection of $a$. Indeed, it is
easy to see that the absolute value of $\f$ is also a supporting
functional of $a$, which must coincide with $\f$ by uniqueness. In
the same way, $e\f e=\f$ for $e\f e$ is again a supporting
functional for $a$. (In fact, one can easily show that $\f$ is
affiliated with the von Neumann subalgebra generated by the
spectral projections of $a$.)

Next, let $E_{(2)}$ be the 2-concavification of $E$. $E_{(2)}$ is
$p/2$-convex and $q/2$-concave (and so  $E_{(2)}(\M)$ is also
uniformly smooth). Consider the one dimensional subspace $\com
a^2\subset  E_{(2)}(\M)$ generated by $a^2$, and the functional
$f_0:\; \com a^2\to\com$ defined by $f_0(\l a^2)=\l$. Then
$\|f_0\|=1$ and $f_0(a^2)=1$. By the Hahn-Banach theorem, $f_0$
extends to a norm one functional $f$ on $E_{(2)}(\M)$. Then $f$ is
the unique  supporting functional of $a^2$, and the preceding
argument shows that $efe=f\ge0$. Let $\psi=af$. We claim that
$\psi$ is a norm one functional on $E(\M)$ and supports $a$.
Indeed, for any $x\in E(\M)$, by the Cauchy-Schwarz inequality
 \be
 |\psi(x)|&=&|\t(xaf)|=\big|\t(f^{1/2}xaf^{1/2})\big|\\
 &\le&
 \big\|f^{1/2}x\big\|_2\,\big\|af^{1/2}\big\|_2=
 \big(f(|x^*|^2)\big)^{1/2}\,\big(f(a^2)\big)^{1/2}\\
 &\le& \big\|\,|x^*|^2\big\|_{E_{(2)}(\M)}^{1/2}=\|x\|_{E(\M)}\,.
 \ee
Thus $\|\psi\|\le1$. However, $\psi(a)=f(a^2)=1$. Then our claim
follows.  Therefore, by uniqueness, $\f=\psi$, i.e. $\f=af$.
Passing to adjoints, we also have $\f=fa$.

Now since  $fa=af$, inequality (\ref{gro thm1}) becomes
 $$\|T(ha+ah)\|^2\le2\f(h^2a + hah)
 =2f(h^2a^2 + haha).$$
On the other hand (recalling that $f\ge0$),
 $$f\big((ha+ah)^2\big)=2f(haha)+ f(h^2a^2) +f(ha^2h)
 \ge f(h^2a^2 + haha).$$
Therefore,
 \beq\label{gro thm2}
 \|T(ha+ah)\|^2\le2f\big((ha+ah)^2\big).
 \eeq
We will apply Lemma \ref{density}. To this end we need to deal
with operators supported by  $e^\perp$.  We claim that $T(x)=0$
for every $x\in E(M)$ such that $e^\perp xe^\perp=x$. It suffices
to consider the case where $x$ is selfadjoint. Then
 \be
 \|T(a)\|^2 + t^2 \|T(x)\|^2
 =\E\|T(a)+ t\e T(x)\|^2
 \le\E\|a+ t\e  x\|^2\,.
 \ee
Since $a$ and $x$ are of disjoint support, by considering the
commutative von Neumann subalgebra generated by $a$ and $x$, we
can assume that $a$ and $x$ are functions of disjoint support.
Thus  the $p$-convexity of $E$ implies that
 \be
 \|a+ t\e  x\|=\|(|a|^p+ t^p |x|^p)^{1/p}\|
 \le(\|a\|^p+ t^p \|x\|^p)^{1/p}\,.
 \ee
Combining the preceding inequalities (recalling that
$\|a\|=\|T(a)\|=1$), we get
 $$t^2 \|T(x)\|^2\le {\rm O}(t^p);$$
whence the claim as $t\to0$ by the assumption that $p>2$.

Now let $h\in \M_h$ and $x=ha+ah+e^\perp he^\perp$. Using the
previous claim, (\ref{gro thm2})  and the fact that $f$ is
supported by $e$, we have
 $$\|T(x)\|^2=\|T(ha+ah)\|^2
 \le 2f\big((ha+ah)^2\big)=2f(x^2)\,.$$
By the density of $\{ha+ah+e^\perp he^\perp\;:\; h\in \M_h\}$ in
$E(\M)_h$ given by Lemma \ref{density}, we deduce that
$\|T(x)\|^2\le 2f(x^2)$ for any selfadjoint $x\in E(\M)$. It then
follows that
 $$\forall\; x\in E(\M)\quad
 \|T(x)\|^2\le 2f(x^*x+xx^*)\,;$$
Namely, (G) holds in this subcase with $K=\sqrt2$.

\medskip\n{\bf Subcase 2: $a\not\ge0$.}
  Let $a=u|a|$ be the polar
decomposition of $a$. Let $e_1=u^*u$ and $e_2=uu^*$. Then $e_1$
and $e_2$ are two equivalent projections of $\M$. Since $\M$ is
finite, their complementary projections $e_1^\perp$ and
$e_2^\perp$ are also equivalent (see \cite[Proposition
V.1.38]{tak-I}).  Therefore, there exists a partial isometry $v\in
\M$ such that $v^*v=e_1^\perp$ and $vv^*=e_2^\perp$. Set $w=u+v$.
Then $w$ is a unitary  and $a=w|a|$.

Now consider a new map $S:\;E(\M)\to \H$ defined by $S(y)= T(wy)$.
Then $S$ has norm one and attains its norm at $|a|$. Therefore, by
Subcase 1.1, there exists a norm one positive functional $g\in
E_{(2)}(\M)^*$ such that
 $$\forall\; y\in E(\M)\quad \|S(y)\|\le 2 g(y^*y+ yy^*);$$
whence (by writing $y=w^*x$)
 $$\forall\; x\in E(\M)\quad \|T(x)\|\le 2 g(x^*x+ w^*xx^*w)
 \le 4f(x^*x+xx^*),$$
where $f=(g+wg w^*)/2$. Therefore, we still have the Grothendieck
factorization for $E(M)$ with $K=2$. Thus the proof of Case 1 is
complete.

\medskip

\n {\bf Case 2:} {\it $E$ is  $p$-convex with constant $1$ for
some $p>2$.} We will show the noncommutative Khintchine inequality
for $F(\M)$, where $F=E'$. To this end note that $F$ is
$p'$-concave with constant 1, where $p'$ denotes the conjugate
index of $p$. In particular, $F$ is order continuous.
Consequently, $\M$ is dense in $F(\M)$ (see \cite[Lemma
4.5]{xu-camb}). Thus in order to prove inequality (K) we need only
to consider finite sequences $(x_k)$ in $\M$.

Now let $r>1$ and consider the $r$-convexification $F^{(r)}$ of
$F$. The order continuity of $F^{(r)}$ implies that for any $x\in
\M$
 $$\lim_{r\to1}\big\|x\big\|_{F^{(r)}(\M)}=\big\|x\big\|_{F(\M)}.$$
Thus we are reduced to show inequality (K) with $F^{(r)}(\M)$ in
place of $F(\M)$ for all $r$ close to 1. However, by Theorem
\ref{equivalence}, this is equivalent to the validity of the
little Grothendieck theorem for $G(M)$, where $G$ is the dual
space of $F^{(r)}$. Since $F^{(r)}$ is $r$-convex and
$rp'$-concave with constant 1, $G$ is $r'$-concave and $s$-convex
with constant 1 , where $s=rp/(1+(r-1)p)$. For $r>1$ sufficiently
close to 1 we still have $s>2$. Thus $G$ verifies the condition of
Case 1. Consequently,  the little Grothendieck theorem holds for
$G(\M)$, so Case 2 is done.

\medskip

\n {\bf Case 3:} {\it The general case.}
 Recall that
$E$ is a 2-convex r.i. space on $[0,1]$.  We will show the
2-concavity inequality (C). To this end fix a map $T:\; E(\M)\to
\H$ with $\|T\|\le 1$.  Let $r>1$ and consider the
$r$-convexification $E^{(r)}$ of $E$. By the H\"older inequality,
$E^{(r)}\subset E$ and the inclusion has norm 1; so
$E^{(r)}(\M)\subset E(\M)$ is also a norm one inclusion. Let $\wt
T=T\circ\iota$, where $\iota$ is the natural inclusion from
$E^{(r)}(\M)$ into $E(\M)$. Thus $\wt T:\;E^{(r)}(\M)\to \H$ is a
contraction. Now $E^{(r)}$ is $2r$-convex with $2r>2$. Therefore,
applying Case 2 to $E^{(r)}(\M)$, we obtain that for any finite
sequence $(x_k)\subset \M$
 $$\big(\sum_k\|\wt T(x_k)\|^2\big)^{1/2}\le K \big\|\big(
 \sum_k x_k^*x_k + x_kx_k^*\big)^{1/2}\big\|_{E^{(r)}(\M)}\,.$$
Namely,
 $$\big(\sum_k\|T(x_k)\|^2\big)^{1/2}\le K \big\|\big(
 \sum_k x_k^*x_k + x_kx_k^*\big)^{r/2}\big\|_{E(\M)}^{1/r}\,.$$
As in Case 2 we also have
 $$\forall\; x\in \M,\; x\ge0\quad
 \lim_{r\to1}\big\|x^r\big\|_{E(\M)}=\big\|x\big\|_{E(\M)}\,.$$
This follows from the order continuity or the Fatou property of
$E$.  Therefore,
 $$\big(\sum_k\|T(x_k)\|^2\big)^{1/2}\le K \big\|\big(
 \sum_k x_k^*x_k + x_kx_k^*\big)^{1/2}\big\|_{E(\M)}.$$
That is, inequality (C) holds for all finite sequences
$(x_k)\subset \M$. To pass from $\M$ to $E(\M)$ we use
approximation as usual in such a situation. Indeed, if $E$ is
order continuous, $\M$ is dense in $E(\M)$, so we are done in this
case. Otherwise, $E=F^*$ with $F=E'$. By Lemma
\ref{approximation}, there exists a net $(T_i)$ of w*-continuous
finite rank maps in the unit ball of $B(E(\M),\; H)$ such that
$T_i\to T$ strongly. Since inequality (C) is stable under strong
limit,  we are reduced to prove (C) for each $T_i$. Replacing $T$
by $T_i$ if necessary, we can assume that $T$ itself is
w*-continuous and of finite rank. Now fix a finite sequence
$(x_k)\subset E(\M)$. For each $n\in\nat$ let
$x_{k,n}=x_k\un_{[0,\;n]}(|x_k|)$, where $\un_{[0,\;n]}(x)$
denotes the spectral projection of a positive operator $x$
corresponding to the interval $[0,\;n]$. Then $x_{k,n}\in\M$, so
by the preceding inequality
 $$\big(\sum_k\|T(x_{k,n})\|^2\big)^{1/2}\le K \big\|\big(
 \sum_k x_{k,n}^*x_{k,n} + x_{k,n}x_{k,n}^*\big)^{1/2}\big\|_{E(\M)}.$$
However,  for each $k$, $x_{k,n}\to x_k$ in $E(\M)$ relative to
the w*-topology as $n\to\8$. It follows that $T(x_{k,n})\to
T(x_k)$ in $H$ by virtue of the w*-continuity of $T$. On the other
hand,
 $$\sum_k x_{k,n}^*x_{k,n} + x_{k,n}x_{k,n}^*
 \le \sum_k x_k^*x_k + x_kx_k^*\,.$$
Therefore, we deduce
 $$\big(\sum_k\|T(x_k)\|^2\big)^{1/2}\le K \big\|\big(
 \sum_k x_k^*x_k + x_kx_k^*\big)^{1/2}\big\|_{E(\M)}\,,$$
as desired. Thus the proof of Theorem \ref{gro thm} is
complete.\cqd


\section{Applications to Schur multipliers}


In this section we present some applications of our little
Grothendieck theorem to Schur multipliers. We  characterize the
Schur multipliers from a 2-convex unitary ideal into a 2-concave
one. Now the von Neumann algebra $\M$ is $B(\el_2)$ and the trace
$\t$ is the usual trace $\Tr$. Accordingly, instead of r.i. spaces
on $(0,\;\8)$, we consider r.i. spaces on $\nat$, i.e. symmetric
sequence spaces. Given a symmetric sequence space $E$, we denote
the associated unitary ideal by $S_E$. Namely,
$S_E=E(B(\el_2),\Tr)$ in the previous notation. Note that if
$E=\el_p$, $S_E$ becomes the usual Schatten class $S_p$. In
particular, $S_1$ is the trace class, $S_\8=B(\el_2)$, and
$S_2=\el_2(\nat^2)$ is the Hilbert-Schmidt class. As usual, the
operators in $S_E$ are represented by infinite matrices. Let
$e_{ij}$ be the canonical matrix units. Then every $x\in S_E$ is
given by an infinite matrix $x=(x_{ij})_{i,j\ge0}$, i.e.
 $$x=\sum_{i,j\ge0}x_{ij}\,e_{ij}\,.$$
Equally, $x$ can be also viewed as a function on $\nat^2$. In the
sequel we will not distinguish an infinite matrix and the
corresponding function on $\nat^2$.

Let $E$ and $F$ be two symmetric sequence spaces. Let
$\f=(\f_{ij})$ be an infinite matrix. We call $\f$ a Schur
multiplier from $S_E$ to $S_F$ if the map $M_\f: x\mapsto
\big(\f_{ij}x_{ij}\big)_{i,j\ge0}$ defines a bounded map from
$S_E$ into $S_F$. More generally, if $X$ and $Y$ are two Banach
spaces of complex functions on $\nat^2$, a Schur multiplier from
$X$ into $Y$ is a function $\f$ on $\nat^2$ such that $M_\f$
induces a bounded map from $X$ into $Y$.

Let $1\le p\le\8$. Recall that $E(\el_p)$ is the space of complex
matrices $\f=(\f_{ij})$ such that the sequence $i\mapsto
\|\f_{i\,\cdot}\|_{\el_p}=(\sum_j|\f_{ij}|^p)^{1/p}$ belongs to
$E$ (with the usual convention for $p=\8$). The norm of $E(\el_p)$
is given by
 $$\|\f\|_{E(\el_p)}
 =\big\|\big(\|\f_{i\,\cdot}\|_{\el_p}\big)_{i\ge0}\big\|_E\,.$$
Let $\,^tE(\el_p)=\{\f\;:\; \,^t\f\in E(\el_p)\}$, equipped with
the natural norm, where $^t\f$ is the transpose of $\f$, i.e.
$^t\f_{ij}=\f_{ji}$. Note that $E(\el_p)$ and $\,^tE(\el_p)$ are
K\"othe function spaces on $\nat^2$. If $X$ and $Y$ are two Banach
spaces of functions on $\nat^2$, $X+Y$ and $X\cap Y$ denote their
sum and intersection, respectively. Recall that the norm of $X+Y$
and $X\cap Y$ are given respectively by
 $$\|z\|_{X+Y}=\inf\big\{\|x\|_X+\|y\|_Y\;:\; z=x+y,\; x\in X,\; y\in
 Y\big\}$$
and
 $$\|z\|_{X\cap Y}=\max\big(\|z\|_X,\;\|z\|_Y\big).$$

\begin{lem}\label{SE to S2}
 Let $E$ be a $2$-convex symmetric sequence space with constant
$1$, and let $\f$ be a function on $\nat^2$.  Then the following
assertions are equivalent
 \begin{enumerate}[\rm i)]
 \item $\f$ is a Schur multiplier from $S_E$ to $S_2$;
 \item $\f$ is a Schur multiplier from $E(\el_2)\cap\,^tE(\el_2)$
to $\el_2(\nat^2)$;
 \item $\f\in G(\el_\8)\,+\,^tG(\el_\8)$, where
 $G=((E_{(2)})')^{(2)}$.
 \end{enumerate}
 Moreover,
  $$\big\|M_\f: S_E\to S_2\big\|\approx
  \big\|M_\f: E(\el_2)\cap\,^tE(\el_2)\to \el_2(\nat^2)\big\|
  =\big\|\f\big\|_{G(\el_\8)+\,^tG(\el_\8)}\,,$$
 where the equivalence constants are controlled by a universal
constant.
 \end{lem}

\pf i) $\Rightarrow$ ii). Let $\f$ be a Schur multiplier from
$S_E$ to $S_2$. Let $x$ be a finite matrix. Then by Theorem
\ref{gro thm} and inequality (C) in Theorem \ref{equivalence},
 \be
 \big\|M_\f(x)\big\|_{S_2}
 &=&\big(\sum_{i,j}|\f_{ij}x_{ij}|^2\big)^{1/2}\\
 &\le& K\|M_\f\|\,\big\|\big[\sum_{i,j}|x_{ij}|^2(e_{ij}^*e_{ij}
 +e_{ij}e_{ij}^*)\big]^{1/2}\big\|_{S_E}\\
 &=& K\|M_\f\|\,\big\|\big[\sum_{i,j}|x_{ij}|^2(e_{jj}
 +e_{ii})\big]^{1/2}\big\|_{S_E}\\
 &\le& 2K\|M_\f\|\,\big\|x\big\|_{E(\el_2)\cap\,^tE(\el_2)}\,.
 \ee
Therefore, $\f$ is a Schur multiplier from
$E(\el_2)\cap\,^tE(\el_2)$ to $\el_2(\nat^2)$.

ii) $\Rightarrow$ i). First observe that $S_E$ embeds
contractively into $E(\el_2)\cap\,^tE(\el_2)$. Indeed, let $x\in
S_E$, and let $a_i=\sum_jx_{ij}e_{ij}$. Then by Theorem
\ref{khintchine}, ii)
 \be
 \|x\|_{E(\el_2)}
 =\big\|\big(\sum_ia_ia_i^*\big)^{1/2}\big\|_{S_E}
 \le\big(\E \big\|\sum_i\e_i\,a_i\big\|_{S(E)}^2\big)^{1/2}
 =\|x\|_{S_E}\,;
 \ee
whence the observation. It then follows that
 $$\big\|M_\f: S_E\to S_2\big\|\le
  \big\|M_\f: E(\el_2)\cap\,^tE(\el_2)\to \el_2(\nat^2)\big\|.$$

ii) $\Leftrightarrow$ iii). Let $X$ be a 2-convex K\"othe function
space on $\nat^2$. Then it is clear that $\f$ is a Schur
multiplier from $X$ to $\el_2(\nat^2)$ iff $\f\in
((X_{(2)})')^{(2)}$. Therefore, the equivalence ii)
$\Leftrightarrow$ iii) follows.\cqd

\begin{rk}\label{E to l2}
 It is clear that the symmetric sequence space $G$
in Lemma \ref{SE to S2}, iii) is equal to the space of multipliers
from $E$ to $\el_2$.
 \end{rk}

Let $G$ and $H$ be two 2-convex symmetric sequence spaces with
constant 1. Define $GH$ by
 $$GH=\{xy\;:\; x\in G,\; y\in H\} \quad\mbox{and}\quad
 \|z\|_{GH}=\inf\{\|x\|_G\,\|y\|_H\;:\; x\in G,\; y\in H\}.$$
It is easy to see that $GH$ is again a symmetric sequence space.

\begin{thm}\label{SE to SF}
 Let $E$ and $F$ be two symmetric sequence spaces. Assume that $E$
and $F$ are respectively $2$-convex and $2$-concave with constant
$1$. Then a function $\f$ on $\nat^2$ is a Schur multiplier from
$S_E$ to $S_F$ iff $\f\in L(\el_\8)\,+\,^tL(\el_\8)$, where $L=GH$
with $G=((E_{(2)})')^{(2)}$ and $H=(((F')_{(2)})')^{(2)}$.
Moreover,
  $$\big\|M_\f: S_E\to S_F\big\|\approx
  \big\|\f\big\|_{L(\el_\8)\,+\,^tL(\el_\8)}$$
 with universal constants.
 \end{thm}

\pf Let $\f=\psi\omega$ with $\psi\in G(\el_\8)$ and $\omega\in
H(\el_\8)$. By Lemma \ref{SE to S2}, $\psi$ is a Schur multiplier
from $S_E$ to $S_2$ and $\omega$ a Schur multiplier from $S_{F'}$
to $S_2$. Passing to adjoint, we see that $\omega$ is also a Schur
multiplier from $S_2$ to $S_{F}$. It follows that $\f$ is a Schur
multiplier from $S_E$ to $S_F$. Consequently, every function in
$L(\el_\8)\,+\,^tL(\el_\8)$ is a Schur multiplier from $S_E$ to
$S_F$.

Conversely, let $\f$ be a Schur multiplier from $S_E$ to $S_F$. To
prove that $\f\in L(\el_\8)\,+\,^tL(\el_\8)$, by the Fatou
property of $L(\el_\8)\,+\,^tL(\el_\8)$, we can clearly assume
that $\f$ is a finite matrix, say an $n\times n$ matrix. Thus in
the remainder of the proof all matrices are assumed to be of order
$n$. Since $F$ and $E'$ are 2-concave with constant 1, by
\cite{TJ-uniconv}, $S_{F}$ and $S_{E}^*$ are of cotype 2 with
universal constants. Therefore, by \cite[Theorem 4.1]{pis-fact},
$M_\f: S_E\to S_F$ factors through a Hilbert space $\mathcal H$ as
$M_\f=VU$ with $\|V\|\,\|U\|\le K\|M_\f\|$, where $K$ is a
universal constant. The point now is to show that we may take
$\mathcal H=S_{2}$ and assume that $U$ and $V$ are Schur
multipliers. The following argument is well-known and standard.

By Theorem \ref{gro thm}, there exists a norm one positive
functional $f\in (S_{E_{(2)}})^*$ such that
 $$\forall\; x\in S_E\quad\|U(x)\|^{2}\le K^{2}f(x^*x+xx^*).$$
Hence
 $$\forall\; x\in S_E\quad\|M_\f(x)\|^{2}\le K^{2}f(x^*x+xx^*).$$
Now let $\e=(\e_i)$  be a Rademacher sequence. Let $D_\e$ be the
diagonal matrix whose diagonal entries are the $\e_i$'s. Let $\e'$
be an independent copy of $\e$ and  $D_{\e'}$  the associated
diagonal matrix. Recall that the norm of $S_F$ is unitary
invariant (in fact, what is needed here is the invariance of the
norm by left and right multiplications by unitary diagonal
matrices). Thus by the previous inequality, for any $x\in S_E$ we
have
 \be
 \|M_\f(x)\|^2
 &=&\|D_\e M_\f(x)D_{\e'}\|^2=\|M_\f(D_\e xD_{\e'})\|^2\\
 &\le& K^2f(D_{\e'}x^*xD_{\e'} + D_{\e}xx^*D_\e)\\
 &=&K^2\big[D_{\e'}fD_{\e'}(x^*x) + D_\e fD_{\e}(xx^*)\big].
 \ee
Taking expectation, we deduce that
 $$\|M_\f(x)\|^2\le K^{2}g(x^*x+xx^*),$$
where $g=\E(D_\e fD_{\e})$. Note that $g\in (S_{E_{(2)}})^*$ is a
positive diagonal matrix, so its diagonal sequence belongs to
$(E_{(2)})'$. The preceding inequality can be rewritten as
 $$\|M_\f(x)\|^2\le K^{2}\big(\big\|g^{1/2}x\big\|_{S_2}^2\,+\,
  \big\|xg^{1/2}\big\|_{S_2}^2\big).$$
It follows that there exist two bounded maps $v$ and $v'$ from
$S_2$ to $S_F$  such that
 $$M_\f=vu + v'u'\,,$$
where $u$ and $u'$ are respectively the left and right
multiplications by $g^{1/2}$. Note that $u=M_\psi$ and
$u'=M_{\psi'}$ with $\psi_{ij}= \a_i$ and $\psi_{ij}'=\a_j$, where
$\a_i= g_{ii}^{1/2}$. Using an average argument as above, we can
further assume that $v$ and $v'$ are also given by Schur
multipliers $M_\omega$ and $M_{\omega'}$. Therefore, $\omega$ and
$\omega'$ are Schur multipliers from $S_2$ to $S_F$, and hence
also from $S_{F'}$ to $S_2$. Thus by Lemma \ref{SE to S2}, $\o,
\o'\in F(\el_\8)\,+\,^tF(\el_\8)$.

Now it is easy to show that $\f=\psi \o+\psi'\o'$ belongs to
$L(\el_\8)\,+\,^tL(\el_\8)$. Indeed, let $\o=\d+\g$ with $\d\in
F(\el_\8)$ and $\g\in\, ^tF(\el_\8)$. It is clear that $\psi\d\in
L(\el_\8)$. We next show that $\psi\g\in
L(\el_\8)\,+\,^tL(\el_\8)$. To this end, by permutations of rows
and columns if necessary, we may assume that the sequence $(\a_i)$
and $(\b_j)$ are nonincreasing, where $\b_j=\sup_i|\g_{ij}|$.
Define $\g'_{ij}=\g_{ij}$ if $i\le j$ and $\g'_{ij}=0$ if $i> j$.
Set $\g''=\g-\g'$. Then $\sup_j|\g'_{ij}|\a_i\le \b_i\a_i$ and
$\sup_i|\g''_{ij}|\a_i\le \b_j\a_j$. It follows that $\psi\g'\in
L(\el_\8)$ and $\psi\g''\in \,^tL(\el_\8)$, so $\psi\g\in
L(\el_\8)\,+\,^tL(\el_\8)$. Consequently, $\psi\o\in
L(\el_\8)\,+\,^tL(\el_\8)$. Similarly, $\psi'\o'\in
L(\el_\8)\,+\,^tL(\el_\8)$. Therefore, $\f\in
L(\el_\8)\,+\,^tL(\el_\8)$. Thus the proof of the theorem is
complete.\cqd

\medskip

The preceding theorem extends the characterization of Schur
multipliers from $S_q$ to $S_p$ for $1\le p\le 2\le q\le\8$ in
\cite{xu-gro} (see also \cite{pisshlyak} for the case of $p=1$ and
$q=\8$). If one of $E$ and $F$ is an $\el_p$, the space $L$ in
Theorem \ref{SE to SF} is easy to be determined. For instance, let
us consider the case where $F=\el_p$ with $1\le p\le2$ (and $E$ is
still 2-convex with constant 1). By Remark \ref{E to l2}, $L=GH$
coincides with the space of multipliers from $E$ to $F$. Thus if
$F=\el_p$, this latter space is equal to $((E_{(p)})')^{(p)}$.
Thus we get the following

\begin{cor}
 Let $E$ be a $2$-convex symmetric sequence space with constant
$1$ and $1\le p\le 2$. Then a function $\f$ on $\nat^2$ is a Schur
multiplier from $S_E$ to $S_p$ iff $\f\in
G_1(\el_\8)\,+\,^tG_1(\el_\8)$, where $G_1=((E_{(p)})')^{(p)}$.
 \end{cor}

The previous arguments apply equally to the case where one of the
unitary ideals $S_E$ and $S_F$ is replaced by a K\"othe function
space on $\nat^2$. By symmetry, it suffices to consider the case
where the second ideal $S_F$ is replaced by a 2-concave K\"othe
function space on $\nat^2$.

\begin{thm}\label{SE to X}
 Let $E$ be a $2$-convex symmetric sequence space with
constant $1$, and let $X$ be a $2$-concave K\"othe function space
on $\nat^2$ with constant $1$. Then a function $\f$ on $\nat^2$ is
a Schur multiplier from $S_E$ to $X$ iff $\f\in
[G(\el_\8)\,+\,^tG(\el_\8)] Y$, where $G=((E_{(2)})')^{(2)}$ and
$Y=(((X')_{(2)})')^{(2)}$. Moreover, the relevant constants are
controlled by a universal constant.

In particular, if $X=\el_p(\nat^2)$ with $1\le p\le2$, then $\f$
is a Schur multiplier from $S_E$ to $\el_p(\nat^2)$ iff $\f\in
G_1(\el_q)\,+\,^tG_1(\el_q)$, where $G_1=((E_{(p)})')^{(p)}$ and
$q=2p/(2-p)$.
 \end{thm}

\pf This proof is almost the same as  that of Theorem \ref{SE to
SF}. The only difference is that the space of Schur multipliers
from $S_2$ to $X$ coincides with the space $Y$, that makes simpler
the present proof. We leave the details to the reader. \cqd

\medskip

 The theorem above in the case of $S_E=B(\el_2)$ and
$X=\el_1(\nat^2)$ goes back to \cite[Example b]{lust-co} (see also
\cite[Theorem 4.1]{ps}).


\end{document}